\documentclass [12pt]{article}
\usepackage{amssymb}
\def\lanbox{\hbox{$\, \vrule height 0.25cm width 0.25cm depth 0.01cm \,$}}

\usepackage{color}

\begin{document}

\centerline{\Large\bf Solvability of
some integro-differential equations}

\centerline{\Large\bf with the logarithmic Laplacian}

\bigskip

\centerline{Vitali Vougalter$^{1}$,   Vitaly Volpert$^2$ $^3$}

\centerline{$^{1}$ Department of Mathematics, University of Toronto,
Toronto, Ontario, M5S 2E4, Canada}

\centerline{e-mail: vitali@math.toronto.edu}

\bigskip

\centerline{$^2$ Institute Camille Jordan, UMR 5208 CNRS, University Lyon 1}

\centerline{Villeurbanne, 69622, France}

\centerline{$^3$ Peoples' Friendship University of Russia, 6 Miklukho-Maklaya
St, Moscow, 117198, Russia}

\centerline{e-mail: volpert@math.univ-lyon1.fr}

\bigskip
\bigskip
\bigskip

\noindent {\bf Abstract.}
We address the existence in the sense of sequences of solutions for a
certain integro-differential type problem involving the logarithmic
Laplacian. The argument is based on the fixed point
technique when such equation contains the
operator without the Fredholm property. It is established that, under the
reasonable technical conditions, the convergence in $L^{1}({\mathbb R}^{d})$ of
the integral kernels yields the existence and convergence in
$L^{2}({\mathbb R}^{d})$ of the solutions.

\bigskip
\bigskip

\noindent {\bf Keywords:} solvability conditions, non-Fredholm
operators, logarithmic Laplacian, integral kernel

\noindent {\bf AMS subject classification:} 35P30, \ 45K05, \ 47G20

\bigskip
\bigskip
\bigskip
\bigskip

\section{Introduction}

 \noindent
 We recall that a linear operator $L$ acting from a Banach
 space $E$ into another Banach space $F$ satisfies the Fredholm
 property if its image is closed, the dimension of its kernel and
 the codimension of its image are finite. Consequently, the
 problem $Lu=f$ is solvable if and only if $\phi_i(f)=0$ for a
 finite number of functionals $\phi_i$ from the dual space $F^*$.
 Such properties of Fredholm operators are broadly used in many
 methods of the linear and nonlinear analysis.

 \noindent
 Elliptic equations in bounded domains with a sufficiently smooth
 boundary satisfy the Fredholm property if the ellipticity
 condition, proper ellipticity and Lopatinskii conditions are
 fulfilled (see e.g. \cite{Ag}, \cite{E09}, \cite{LM}, \cite{Volevich}).
 This is the main result of the theory of linear elliptic problems.
 In the situation of unbounded domains, such conditions may
 not be sufficient and the Fredholm property may not be satisfied.
 For instance, the Laplace operator, $Lu = \Delta u$, in $\mathbb R^d$
 fails to satisfy the Fredholm property when considered in
 H\"older spaces, $L : C^{2+\alpha}(\mathbb R^d) \to C^{\alpha}(\mathbb
 R^d)$, or in Sobolev spaces,  $L : H^2(\mathbb R^d) \to L^2(\mathbb
 R^d)$.

 \noindent
 Linear elliptic equations in unbounded domains satisfy the
 Fredholm property if and only if, in addition to the conditions
 given above, the limiting operators are invertible (see \cite{V11}).
 In some simple cases, the limiting operators can be constructed
 explicitly. For example, if
 $$
 L u = a(x) u'' + b(x) u' + c(x) u , \;\;\; x \in \mathbb R ,
 $$
 where the coefficients of the operator have limits at the infinities,

 $$ a_\pm =\lim_{x \to \pm \infty} a(x) , \;\;\;
 b_\pm =\lim_{x \to \pm \infty} b(x) , \;\;\;
 c_\pm =\lim_{x \to \pm \infty} c(x) , $$
 the limiting operators are equal to:

 $$ L_{\pm}u = a_\pm u'' + b_\pm u' + c_\pm u . $$
 Since the coefficients here are constants, the essential spectrum of
 such operator, that is the set of complex numbers $\lambda$ for
 which the operator $L-\lambda$ fails to satisfy the Fredholm
 property, can be explicitly found by means of the Fourier
 transform:
 $$
 \lambda_{\pm}(\xi) = -a_\pm \xi^2 + b_\pm i\xi + c_\pm , \;\;\;
 \xi \in \mathbb R .
 $$
 The invertibility of the limiting operators is equivalent to the condition
 that the essential spectrum does not contain the origin.

 \noindent
 In the case of general elliptic equations, the same assertions are valid.
 The Fredholm property is satisfied if the origin does not belong to the
 the essential spectrum or if the limiting operators are invertible. However,
 these conditions may not be explicitly formulated.

 \noindent
 In the situation of non-Fredholm operators the usual solvability
 relations may not be applicable and solvability conditions
 are, in general, not known. There are some classes of operators
 for which solvability relations are obtained. Let us illustrate
 them with the following example. Consider the equation
 \begin{equation}
 \label{int1}
 Lu \equiv \Delta u + a u = f
 \end{equation}
 in $\mathbb R^d, \ d\in {\mathbb N}$, where $a$ is a positive constant.
 The operator $L$ here coincides with its limiting operators. The
 homogeneous equation has a nontrivial bounded solution. Thus, the
 Fredholm property is not satisfied. However, since the operator
 has constant coefficients, we can apply the Fourier transform and
 find the solution explicitly. The solvability conditions can be
 formulated as follows. If $f \in L^2(\mathbb R^d)$ and 
 $xf \in L^1(\mathbb R^d)$, then there
 exists a unique solution of such problem in $H^2(\mathbb R^d)$ if and
 only if
 $$  
 \Bigg(f(x),\frac{e^{ipx}}{(2\pi)^{\frac{d}{2}}}\Bigg)_{L^2(\mathbb R^d)}=0, \quad 
 p\in S_{\sqrt{a}}^{d} \quad a.e.
 $$
 (see \cite{VV103}). Throughout the article $S_{r}^{d}$ denotes the sphere
 in $\mathbb R^d$ of radius $r$ centered at the origin.
 Thus, though the operator fails to satisfy the Fredholm property,
 the solvability relations are stated similarly. However,
 such similarity is only formal because the range of the operator is
 not closed.

 \noindent
 In the case of the operator involving a scalar potential,
 $$
 L u \equiv \Delta u + a(x) u = f ,
 $$
 the Fourier transform is not directly applicable. Nevertheless, the solvability
 conditions in ${\mathbb R}^{3}$ can be obtained by a rather sophisticated 
 application of the theory of self-adjoint
 operators (see \cite{VV08}). As before, the solvability relations are
 given in terms of the orthogonality to the solutions of the homogeneous
 adjoint equation. There are several other examples of linear elliptic
 non-Fredholm operators for which the solvability
 conditions can be derived (see ~\cite{V11}, ~\cite{VKMP02}, ~\cite{VV131},
 ~\cite{VV103}).

 \noindent
 Solvability relations play a crucial role in the analysis of
 nonlinear elliptic problems. In the case of the operators without the Fredholm
 property, in spite of a certain progress in the understanding of the linear
 equations, there exist only few examples where nonlinear non-Fredholm
 operators are studied (see ~\cite{DMV05}, ~\cite{DMV08}, ~\cite{EV20},
 ~\cite{EV211}, ~\cite{VV11}, ~\cite{VV103}, ~\cite{VV18}, ~\cite{VV21}). The
 work ~\cite{E10} addresses the
 finite and infinite dimensional attractors for the evolution
 equations of mathematical physics. The large time behavior of solutions
 of a class of fourth-order parabolic problems defined on unbounded domains
 using the Kolmogorov $\varepsilon$-entropy as a measure was studied
 in ~\cite{EP07}. In ~\cite{EZ01} the authors discuss the attractor for a
 nonlinear reaction-diffusion system in an unbounded domain in
 ${\mathbb R}^{3}$. The articles ~\cite{GS05} and ~\cite{RS01} are important
 for the understanding of the Fredholm and properness properties of quasilinear
 elliptic systems of second order and of the operators of this kind on
 ${\mathbb R}^{N}$. The work ~\cite{GS10} is dedicated to the exponential
 decay and the Fredholm properties in second order quasilinear elliptic
 systems. The article ~\cite{E18} is a systematic study of a dynamical
 systems approach to investigating the symmetrization and stabilization
 properties
 of the nonnegative solutions of nonlinear elliptic equations in asymptotically
 symmetric unbounded domains.

 \noindent
 In the present work we consider another class of stationary
 nonlinear equations, for which the Fredholm property may not be satisfied:
\begin{equation}
\label{id1}
\Big[-\frac{1}{2}\hbox{ln}(-\Delta)\Big]u+au + 
\int_{{\mathbb R}^{d}}G(x-y)F(u(y),y)dy = 0, \quad x\in {\mathbb R}^{d}, \quad
d\in {\mathbb N},
\end{equation}
where $a\in {\mathbb R}$ is a constant.

\noindent
The logarithmic Laplacian
$\hbox{ln}(-\Delta)$ is the operator with Fourier symbol $2\hbox{ln}|p|$.
It arises as the formal derivative $\partial_{s}|_{s=0}(-\Delta)^{s}$ of
fractional Laplacians at $s=0$. The operator $(-\Delta)^{s}$ is extensively
used, for instance in the studies of
the anomalous diffusion problems (see e.g. ~\cite{MK00}, ~\cite{VV21} and the
references therein). Spectral properties of the logarithmic Laplacian in
an open set of finite measure with Dirichlet boundary conditions were considered
in ~\cite{LW21} (see also ~\cite{CW19}). The studies of $\hbox{ln}(-\Delta)$
are important for the understanding of the asymptotic spectral properties of
the family of fractional Laplacians in the limit $s\to 0^{+}$. In ~\cite{JSW20}
it was established that such operator allows to characterize the $s$-dependence
of solution to fractional Poisson equations for the full range of exponents
$s\in (0, 1)$. The solvability of the integro-differential equation similar
to (\ref{id1}) with the standard Laplace operator was discussed in ~\cite{VV11}.

\noindent
The operator involved in problem (\ref{id1}) is given by
\begin{equation}
\label{la}
l_{a}:=\frac{1}{2}\hbox{ln}(-\Delta)-a, \quad a\in {\mathbb R}  
\end{equation}
and is considered on $L^{2}({\mathbb R}^{d}), \ d\in {\mathbb N}$. By virtue of
the standard Fourier transform, it can be trivially obtained that the
essential spectrum of (\ref{la}) equals to
\begin{equation}
\label{ess}
\lambda_{a}(p)=\hbox{ln}|p|-a, \quad p\in {\mathbb R}^{d}, \quad
a\in {\mathbb R}.
\end{equation}  
As distinct from many previous articles devoted to the operators without the
Fredholm property, (\ref{ess}) fills not a semi-axis but the whole real line,
which can be easily observed from (\ref{ess}). Hence, the inverse of (\ref{la})
is unbounded.

\noindent
The article ~\cite{EV23} deals with the solvability of certain linear
nonhomogeneous equations containing the logarithmic Laplacian and the
logarithmic Schr\"odinger operator with a shallow, short-range scalar
potential. These results were generalized to the case of the logarithmic
Schr\"odinger operators in higher dimensions in ~\cite{EV231}. The work
~\cite{EV24} is devoted to the solvability of linear and nonlinear
nonhomogeneous equations in one dimension involving the logarithmic Laplacian
and the transport term. The logarithmic Schr\"odinger operator and associated
Dirichlet problems were covered in ~\cite{F23}. 
Symmetry of
positive solutions for Lane-Emden systems involving the logarithmic Laplacian
was discussed in ~\cite{ZKR23}. The logarithmic Laplacian is also used in the
geometric context of the $0$-fractional perimeter (see ~\cite{DNP21}).

\noindent
In the population dynamics the 
integro-differential problems describe the models with the intra-specific
competition and the nonlocal consumption of resources (see e.g.
~\cite{ABVV10}, ~\cite{BNPR09}). Solving them is important for the
understanding of the
emergence and propagation of patterns in the theory of speciation
(see ~\cite{VV13}). Standing lattice solitons in the discrete NLS equation
with saturation were covered in ~\cite{AKLP19}.

\noindent
Let us use the explicit form of
the solvability conditions and explore the existence of solutions of our
nonlinear equation.


\setcounter{equation}{0}

\section{Formulation of the results}

The nonlinear part of our problem (\ref{id1}) will satisfy the
following regularity conditions.

\bigskip

\noindent
{\bf Assumption 1.} {\it Let $d\in {\mathbb N}$. Function
$F(u,x): {\mathbb R}\times {\mathbb R}^{d}\to {\mathbb R}$ is satisfying the
Caratheodory condition (see ~\cite{K64}), such that
\begin{equation}
\label{ub1}
|F(u,x)|\leq k|u|+h(x) \quad for \quad u\in {\mathbb R}, \quad
x\in {\mathbb R}^{d}
\end{equation}
with a constant $k>0$ and $h(x):{\mathbb R}^{d}\to {\mathbb R}^{+}, \quad
h(x)\in L^{2}({\mathbb R}^{d})$. Moreover, it is a Lipschitz continuous function,
so that
\begin{equation}
\label{lk1}
|F(u_{1},x)-F(u_{2},x)|\leq l |u_{1}-u_{2}| \quad for \quad any \quad
 u_{1,2}\in{\mathbb R}, \quad x\in {\mathbb R}^{d}
\end{equation}
with a constant $l>0$.}

\bigskip

\noindent
The work ~\cite{BO86} deals with the solvability of a local elliptic
equation in a bounded domain in ${\mathbb R}^{N}$. The nonlinear function
involved there was allowed to have a sublinear growth.

\noindent
In order to establish the existence of the solution of problem (\ref{id1}), we
will use the auxiliary equation
\begin{equation}
\label{ae1}
\Big[\frac{1}{2}\hbox{ln}(-\Delta)\Big]u-au=\int_{{\mathbb R}^{d}}G(x-y)
F(v(y),y)dy, \quad x\in {\mathbb R}^{d}, \quad d\in {\mathbb N},
\end{equation}
where $a\in {\mathbb R}$ is a constant.

\noindent
Let us define
\begin{equation}
\label{ip}  
(f_{1}(x),f_{2}(x))_{L^{2}({\mathbb R}^{d})}:=\int_{{\mathbb R}^{d}}f_{1}(x)\bar{f_{2}}(x)dx,
\end{equation}
with a slight abuse of notations when such functions do not belong to
$L^{2}({\mathbb R}^{d})$, like for instance those contained in orthogonality
relations (\ref{or1}) and (\ref{or2}) below. Clearly, if
$f_{1}(x)\in L^{1}({\mathbb R}^{d})$ and
$f_{2}(x)\in L^{\infty}({\mathbb R}^{d})$, then the integral in the right side of
(\ref{ip}) is well defined.

\noindent
The main issue for our equation (\ref{ae1}) is that 
we deal with the self-adjoint, non-Fredholm  operator (\ref{la}) involved in it,
which is the obstacle to solve our problem.
The similar situations but
in linear problems, both self-adjoint and non self-adjoint
containing the differential operators without the Fredholm property
were considered extensively in recent years (see
~\cite{V11}, ~\cite{VKMP02}, ~\cite{VV08}, ~\cite{VV103}).
 
\noindent
We establish that under the reasonable auxiliary conditions,
problem (\ref{ae1}) defines a map $t_{a}:
L^{2}({\mathbb R}^{d})\to L^{2}({\mathbb R}^{d})$ with the constant
$a\in {\mathbb R}$, which is a strict contraction.

\bigskip

\noindent
{\bf Theorem 1.}  {\it Let $G(x): {\mathbb R}^{d}\to
{\mathbb R}, \ G(x)\in L^{1}({\mathbb R}^{d}), \ xG(x)\in L^{1}({\mathbb R}^{d})$
with $d\in {\mathbb N}, \ a\in {\mathbb R}$ and Assumption 1 is valid.

\medskip

\noindent
We also assume that orthogonality conditions (\ref{or1}) hold when $d=1$
and (\ref{or2}) are valid for $d\geq 2$ and
$(2 \pi)^{\frac{d}{2}}N_{a}l<1$ with $N_{a}$ defined in (\ref{Na}).
Then the map $v \mapsto t_{a}v=u$ on $L^{2}({\mathbb R}^{d})$ defined by equation
(\ref{ae1}) has a unique fixed point $v_{a}$, which is the only
solution of problem (\ref{id1}) in $L^{2}({\mathbb R}^{d})$.

\medskip

\noindent
The fixed point
$v_{a}, \ a\in {\mathbb R}$ is nontrivial in the whole space provided the
intersection of supports of the Fourier transforms
of functions
$$
supp\widehat{F(0,x)}\cap supp \widehat{G}
$$
is a set of nonzero
Lebesgue measure in ${\mathbb R}^{d}$.}

\bigskip

\noindent
Let us introduce the sequence of
approximate equations related to problem (\ref{id1}), such that
\begin{equation}
\label{id1m}
\Big[-\frac{1}{2}\hbox{ln}(-\Delta)\Big]u_{m}+au_{m} + 
\int_{{\mathbb R}^{d}}G_{m}(x-y)F(u_{m}(y),y)dy = 0, \quad x\in {\mathbb R}^{d},
\quad d\in {\mathbb N}
\end{equation}
with a constant $a\in {\mathbb R}$ and $m\in {\mathbb N}$.
The sequence of kernels $\displaystyle{\{G_{m}(x)\}_{m=1}^{\infty}}$ tends
to $G(x)$ as $m\to \infty$ in the appropriate function spaces specified below.
We establish that, under the certain technical conditions, each of
equations (\ref{id1m}) admits a unique solution
$u_{m}(x)\in L^{2}({\mathbb R}^{d})$, the limiting problem (\ref{id1})
possesses a unique solution $u(x)\in L^{2}({\mathbb R}^{d})$, and
$u_{m}(x)\to u(x)$ in $L^{2}({\mathbb R}^{d})$ as $m\to \infty$. This is the
so-called {\it existence of solutions in the sense of sequences}
(see ~\cite{VV18}, also ~\cite{V11}). In this
case, the solvability conditions can be stated for the iterated kernels
$G_{m}$. They imply the convergence of the kernels in terms of the Fourier
transforms (see the Appendix) and, consequently, the convergence of the
solutions (Theorem 2).
The similar ideas in the sense of the standard Schr\"odinger type operators
were used in ~\cite{VV131}. Our second main statement is as
follows.

\bigskip

\noindent
{\bf Theorem 2.}  {\it Let $m\in {\mathbb N}, \
G_{m}(x): {\mathbb R}^{d}\to {\mathbb R}, \ G_{m}(x)\in L^{1}({\mathbb R}^{d}), \
xG_{m}(x)\in L^{1}({\mathbb R}^{d})$ with $d\in {\mathbb N}, \ a\in {\mathbb R}$,
such that
$G_{m}(x)\to G(x)$ in $L^{1}({\mathbb R}^{d})$ and
$xG_{m}(x)\to xG(x)$ in $L^{1}({\mathbb R}^{d})$
as $m\to \infty$ and Assumption 1  holds.

\medskip

\noindent
Let us also assume that orthogonality relations (\ref{or3}) are valid
for $d=1$ and orthogonality conditions (\ref{or4}) hold when $d\geq 2$ and
$$ 
(2\pi)^{\frac{d}{2}}N_{a, \ m}l\leq 1-\varepsilon
$$
for all $m\in {\mathbb N}$ with a certain fixed $0<\varepsilon<1$, where
$N_{a, \ m}$ is introduced in (\ref{Nam}).

\medskip

\noindent
Then each equation (\ref{id1m}) has a unique solution
$u_{m}(x)\in L^{2}({\mathbb R}^{d})$, and limiting problem (\ref{id1})
admits a unique solution $u(x)\in L^{2}({\mathbb R}^{d})$.

\medskip

\noindent
Moreover, $u_{m}(x)\to u(x)$ in $L^{2}({\mathbb R}^{d})$
as $m\to \infty$.

\medskip

\noindent
The unique solution $u_{m}(x)$ of each equation (\ref{id1m}) does not
vanish identically in the whole space
provided that the intersection of supports of the Fourier transforms
of functions
$$
supp\widehat{F(0,x)}\cap supp \widehat{G_{m}}
$$
is a set of nonzero
Lebesgue measure in ${\mathbb R}^{d}$. Analogously, the unique solution $u(x)$
of limiting problem (\ref{id1}) is nontrivial if
$$
supp\widehat{F(0,x)}\cap supp \widehat{G}
$$
is a set of nonzero Lebesgue
measure in ${\mathbb R}^{d}$.}

\bigskip

\noindent
{\bf Remark 1.} {\it The importance of Theorem 2 above is the continuous
dependence of the solutions with respect to the integral kernels.}      


\setcounter{equation}{0}

\section{Solvability of the integro-differential equation}

\bigskip

{\it Proof of Theorem 1.} Let us first suppose that 
for some $v(x)\in L^{2}({\mathbb R}^{d})$ there
exist two solutions $u_{1,2}(x)\in L^{2}({\mathbb R}^{d})$ of equation
(\ref{ae1}). Then the difference function
$w(x):=u_{1}(x)-u_{2}(x)\in  L^{2}({\mathbb R}^{d})$ satisfies the homogeneous
equation
\begin{equation}
\label{hom}  
\Big[\frac{1}{2}\hbox{ln}(-\Delta)\Big]w-aw=0.
\end{equation}
Since the operator in the left side of (\ref{hom}) introduced in (\ref{la})
does not have any nontrivial square integrable zero modes,
$w(x)$ vanishes in  ${\mathbb R}^{d}$.

\noindent
We choose arbitrarily $v(x)\in L^{2}({\mathbb R}^{d})$ and apply the standard
Fourier transform (\ref{ft}) to both sides of problem (\ref{ae1}). This gives
us
\begin{equation}
\label{f1}
\widehat{u}(p)=(2\pi)^{\frac{d}{2}}\frac{\widehat{G}(p)\widehat{\varphi}(p)}
{\hbox{ln}\Big(\frac{|p|}{e^{a}}\Big)}, 
\end{equation}
where $\widehat{\varphi}(p)$ stands for the Fourier image of $F(v(x),x)$.

\noindent
Evidently, the upper bound
$$
|\widehat{u}(p)|\leq (2\pi)^{\frac{d}{2}}N_{a}|\widehat{\varphi}(p)|  
$$
holds.
Clearly, $N_{a}<\infty$ by virtue of Lemma A1 of the Appendix 
under orthogonality relations (\ref{or1}) when $d=1$ and under (\ref{or2})
for $d\geq 2$.

\noindent
This allows us to derive the estimate from above on the norm as
\begin{equation}
\label{l2nub}  
\|u\|_{L^{2}({\mathbb R}^{d})}\leq
(2\pi)^{\frac{d}{2}} N_{a}\|F(v(x),x)\|_{L^{2}({\mathbb R}^{d})}.
\end{equation}
Let us recall inequality (\ref{ub1}) of Assumption 1 above.
Hence, under the given conditions the right side of bound (\ref{l2nub}) is
finite.
This means that for an arbitrary
$v(x)\in L^{2}({\mathbb R}^{d})$ there exists a unique solution
$u(x)\in L^{2}({\mathbb R}^{d})$ of equation (\ref{ae1}) with its Fourier
transform given by (\ref{f1}), such that the map
$t_{a}:  L^{2}({\mathbb R}^{d})\to  L^{2}({\mathbb R}^{d})$ is well defined.

\noindent
This enables us to choose arbitrary $v_{1,2}(x)\in L^{2}({\mathbb R}^{d})$,
so that their images $u_{1,2}:=t_{a}v_{1,2}\in L^{2}({\mathbb R}^{d})$.
From (\ref{ae1}) we easily deduce that
\begin{equation}
\label{aer1}
\Big[\frac{1}{2}\hbox{ln}(-\Delta)\Big]u_{1}-au_{1}=\int_{{\mathbb R}^{d}}
G(x-y)F(v_{1}(y),y)dy,
\end{equation}
\begin{equation}
\label{aer2}
\Big[\frac{1}{2}\hbox{ln}(-\Delta)\Big]u_{2}-au_{2}=\int_{{\mathbb R}^{d}}
G(x-y)F(v_{2}(y),y)dy.
\end{equation}
Let us apply the standard Fourier transform (\ref{ft}) to both sides of
equations (\ref{aer1}) and (\ref{aer2}). This yields
$$
\widehat{u}_{1}(p)=(2\pi)^{\frac{d}{2}}\frac{\widehat{G}(p)\widehat{\varphi_{1}}(p)}
{\hbox{ln}\Big(\frac{|p|}{e^{a}}\Big)}, \quad 
\widehat{u}_{2}(p)=(2\pi)^{\frac{d}{2}}\frac{\widehat{G}(p)\widehat{\varphi_{2}}(p)}
{\hbox{ln}\Big(\frac{|p|}{e^{a}}\Big)}, 
$$
where $\widehat{\varphi_{1}}(p)$ and  $\widehat{\varphi_{2}}(p)$ denote the
Fourier transforms of $F(v_{1}(x),x)$  and $F(v_{2}(x),x)$ respectively.

\noindent
It can be trivially checked that the estimates from above
$$
|\widehat{u}_{1}(p)-\widehat{u}_{2}(p)|\leq (2\pi)^{\frac{d}{2}}N_{a}
|\widehat{\varphi_{1}}(p)-\widehat{\varphi_{2}}(p)|
$$
is valid, which allows us to obtain the upper bound for the norm as
$$
\|u_{1}-u_{2}\|_{L^{2}({\mathbb R}^{d})}\leq (2\pi)^{\frac{d}{2}} N_{a}
\|F(v_{1}(x),x)-F(v_{2}(x),x)\|_{L^{2}({\mathbb R}^{d})}.
$$
We recall inequality (\ref{lk1}) of Assumption 1. Thus,
$$
\|F(v_{1}(x),x)-F(v_{2}(x),x)\|_{L^{2}({\mathbb R}^{d})}\leq l
\|v_{1}(x)-v_{2}(x)\|_{L^{2}({\mathbb R}^{d})},
$$
so that 
\begin{equation}
\label{Tab}  
\|t_{a}v_{1}-t_{a}v_{2}\|_{L^{2}({\mathbb R}^{d})}\leq (2\pi)^{\frac{d}{2}}
N_{a}l\|v_{1}-v_{2}\|_{L^{2}({\mathbb R}^{d})}.
\end{equation}
The constant in the right side of (\ref{Tab}) is less than one as we
assume. By virtue of the Fixed Point Theorem,
there exists a unique function
$v_{a}\in L^{2}({\mathbb R}^{d})$ with the property $t_{a}v_{a}=v_{a}$.
This is the only solution of equation (\ref{id1}) in
$L^{2}({\mathbb R}^{d})$.

\noindent
Let us suppose $v_{a}(x)$ vanishes identically in ${\mathbb R}^{d}$.
This will contradict to our assumption that the Fourier
transforms of $G(x)$ and $F(0,x)$ are nontrivial on a set of nonzero Lebesgue
measure in ${\mathbb R}^{d}$.
\hfill\lanbox

\bigskip

\noindent
Then we turn our attention to establishing the existence in the sense of
sequences of the solutions for our integro-differential equation in
${\mathbb R}^{d}$.

\bigskip

\noindent
{\it Proof of Theorem 2.} According to the result of Theorem 1 above, each
problem (\ref{id1m}) admits a unique solution
$u_{m}(x)\in L^{2}({\mathbb R}^{d}), \ m\in {\mathbb N}$. Limiting equation
(\ref{id1}) has a unique solution $u(x)\in L^{2}({\mathbb R}^{d})$ by virtue
of Lemma A2 below along with Theorem 1.

\noindent
We apply the standard
Fourier transform (\ref{ft}) to both sides of (\ref{id1}) and (\ref{id1m})
and arrive at
\begin{equation}
\label{ump}
\widehat{u}(p)=(2\pi)^{\frac{d}{2}}\frac{\widehat{G}(p)\widehat{f}(p)}
{\hbox{ln}\Bigg(\frac{|p|}{e^{a}}\Bigg)},
\quad \widehat{u_{m}}(p)=(2\pi)^{\frac{d}{2}}\frac{\widehat{G_{m}}(p)
\widehat{f_{m}}(p)}{\hbox{ln}\Bigg(\frac{|p|}{e^{a}}\Bigg)}, \quad
m\in {\mathbb N},
\end{equation}  
where $\widehat{f}(p)$ and $\widehat{f_{m}}(p)$ denote the
Fourier images of $F(u(x), x)$ and $F(u_{m}(x), x)$ respectively.

\noindent
Evidently,
$$
|\widehat{u_{m}}(p)-\widehat{u}(p)|\leq (2\pi)^{\frac{d}{2}}\Bigg\|
\frac{\widehat{G_{m}}(p)}{\hbox{ln}\Bigg(\frac{|p|}{e^{a}}\Bigg)}-
\frac{\widehat{G}(p)}{\hbox{ln}\Bigg(\frac{|p|}{e^{a}}\Bigg)}\Bigg\|_
{L^{\infty}({\mathbb R}^{d})}|\widehat{f}(p)|+
$$
$$
(2\pi)^{\frac{d}{2}}N_{a, \ m}|\widehat{f_{m}}(p)-\widehat{f}(p)|
$$
with $N_{a, \ m}$ defined in (\ref{Nam}). Then
$$ 
\|u_{m}-u\|_{L^{2}({\mathbb R}^{d})}\leq (2\pi)^{\frac{d}{2}}
\Bigg\|\frac{\widehat{G_{m}}(p)}{\hbox{ln}\Bigg(\frac{|p|}{e^{a}}\Bigg)}-
\frac{\widehat{G}(p)}{\hbox{ln}\Bigg(\frac{|p|}{e^{a}}\Bigg)}\Bigg\|_{L^{\infty}
({\mathbb R}^{d})}\|F(u(x), x)\|_{L^{2}({\mathbb R}^{d})}+
$$
\begin{equation}
\label{umul2nub} 
(2\pi)^{\frac{d}{2}}N_{a, \ m}
\|F(u_{m}(x), x)-F(u(x), x)\|_{L^{2}({\mathbb R}^{d})}.
\end{equation}
Let us recall bound (\ref{lk1}) of Assumption 1 above. Hence,
\begin{equation}
\label{Fumu}  
\|F(u_{m}(x), x)-F(u(x), x)\|_{L^{2}({\mathbb R}^{d})}\leq l\|u_{m}(x)-u(x)\|_{L^{2}
({\mathbb R}^{d})}.
\end{equation}
By means of (\ref{umul2nub}) and (\ref{Fumu}), we derive 
$$
\|u_{m}(x)-u(x)\|_{L^{2}({\mathbb R}^{d})}\{1-(2\pi)^{\frac{d}{2}}N_{a, \ m}l\}
\leq
$$
\begin{equation}
\label{umul2nbr}  
(2\pi)^{\frac{d}{2}}
\Bigg\|\frac{\widehat{G_{m}}(p)}{\hbox{ln}\Big(\frac{|p|}{e^{a}}\Big)}-
\frac{\widehat{G}(p)}
{\hbox{ln}\Big(\frac{|p|}{e^{a}}\Big)}\Bigg\|_{L^{\infty}({\mathbb R}^{d})}
\|F(u(x), x)\|_{L^{2}({\mathbb R}^{d})}.
\end{equation}
Using (\ref{umul2nbr}) and (\ref{2rnabml0}), we have
$$
\|u_{m}(x)-u(x)\|_{L^{2}({\mathbb R}^{d})}\leq \frac{(2\pi)^{\frac{d}{2}}}{\varepsilon}
\Bigg\|\frac{\widehat{G_{m}}(p)}{\hbox{ln}\Big(\frac{|p|}{e^{a}}\Big)}-
\frac{\widehat{G}(p)}
{\hbox{ln}\Big(\frac{|p|}{e^{a}}\Big)}\Bigg\|_{L^{\infty}({\mathbb R}^{d})}
\|F(u(x), x)\|_{L^{2}({\mathbb R}^{d})}.
$$
By virtue of inequality (\ref{ub1}) of our Assumption 1, if
$u(x)\in L^{2}({\mathbb R}^{d})$ then
$F(u(x), x)\in L^{2}({\mathbb R}^{d})$ as well. Therefore,
\begin{equation}
\label{umuc}
u_{m}(x)\to u(x), \quad m\to \infty   
\end{equation}
in $L^{2}({\mathbb R}^{d})$ via the result of Lemma A2 of the Appendix.

\noindent
We suppose that the unique solution $u_{m}(x)$ of problem (\ref{id1m})
considered above is trivial in ${\mathbb R}^{d}$ for some $m\in {\mathbb N}$.
This will contradict to our condition that the Fourier transforms of
$G_{m}(x)$ and $F(0, x)$ do not vanish identically 
on a set of nonzero Lebesgue measure in ${\mathbb R}^{d}$. The similar reasoning
holds for the unique solution $u(x)$ of limiting equation (\ref{id1}).
\hfill\lanbox

\bigskip


\setcounter{equation}{0}

\section{Appendix}

\bigskip

Let $G(x)$ be a function,
$G(x): {\mathbb R}^{d}\to {\mathbb R}, \ d\in {\mathbb N}$.
We denote its standard Fourier transform using the hat symbol as
\begin{equation}
\label{ft}  
\widehat{G}(p):={1\over (2\pi)^{\frac{d}{2}}}\int_{{\mathbb R}^{d}}G(x)
e^{-ipx}dx, \quad p\in {\mathbb R}^{d},
\end{equation}
such that
\begin{equation}
\label{inf1}
\|\widehat{G}(p)\|_{L^{\infty}({\mathbb R}^{d})}\leq {1\over (2\pi)^{\frac{d}{2}}}
\|G(x)\|_{L^{1}({\mathbb R}^{d})}
\end{equation}
and
$\displaystyle{G(x)={1\over (2\pi)^{\frac{d}{2}}}\int_{{\mathbb R}^{d}}
\widehat{G}(q)e^{iqx}dq, \ x\in {\mathbb R}^{d}.}$

\noindent
For the technical purposes we introduce the auxiliary quantity
\begin{equation}
\label{Na}
N_{a}:=
\Bigg\|\frac{\widehat{G}(p)}{\hbox{ln}\Big(\frac{|p|}{e^{a}}\Big)}\Bigg\|_
{L^{\infty}({\mathbb R}^{d})},
\end{equation}
where $a\in {\mathbb R}$ is a constant.

\bigskip

\noindent
{\bf Lemma A1.} {\it Let $G(x): {\mathbb R}^{d}\to {\mathbb R}, \ G(x)\in
L^{1}({\mathbb R}^{d}), \ xG(x)\in L^{1}({\mathbb R}^{d})$ with
$d\in {\mathbb N}$ and $a\in {\mathbb R}$.

\medskip
  
\noindent  
a) If $d=1$ then $N_{a}<\infty$ if and only if the orthogonality relations
\begin{equation}
\label{or1}
\Big(G(x), \frac{e^{\pm ie^{a}x}}{\sqrt{2\pi}}\Big)_{L^{2}({\mathbb R})}=0
\end{equation}
are valid.

\medskip

\noindent
b) If $d\geq 2$ then $N_{a}<\infty$ if and only if the orthogonality conditions
\begin{equation}
\label{or2}
\Big(G(x), \frac{e^{ipx}}{(2\pi)^{\frac{d}{2}}}\Big)_{L^{2}({\mathbb R}^{d})}=0, \quad
p\in S_{e^{a}}^{d}
\end{equation}
hold.}

\medskip

\noindent
{\it Proof.} First of all, we discuss the case b) of the
lemma, such that the dimension of the problem $d\geq 2$.

\noindent
Let us introduce the spherical layer for the technical purposes, namely
\begin{equation}
\label{ad}  
A_{\delta}:=\{p\in {\mathbb R}^{d} \ | \ e^{a}(1-\delta)\leq |p|\leq
e^{a}(1+\delta)\}, \quad 0<\delta<1, \quad d\geq 2,
\end{equation}
such that 
\begin{equation}
\label{ghpln}
\frac{\widehat{G}(p)}{\hbox{ln}\Big(\frac{|p|}{e^{a}}\Big)}=
\frac{\widehat{G}(p)}{\hbox{ln}\Big(\frac{|p|}{e^{a}}\Big)}\chi_{A_{\delta}}+
\frac{\widehat{G}(p)}{\hbox{ln}\Big(\frac{|p|}{e^{a}}\Big)}\chi_{A_{\delta}^{c}}.
\end{equation}
Here and below $A^{c}$ will denote the complement of a set
$A\subseteq {\mathbb R}^{d}$. The characteristic function of a set $A$ is being
designated as $\chi_{A}$.

\noindent
Evidently,
$$
A_{\delta}^{c}=A_{\delta}^{c +}\cup A_{\delta}^{c -},
$$
where
\begin{equation}
\label{adcp}  
A_{\delta}^{c +}:=\{p\in {\mathbb R}^{d} \ | \ |p|>e^{a}(1+\delta)\},
\end{equation}
\begin{equation}
\label{adcm}  
A_{\delta}^{c -}:=\{p\in {\mathbb R}^{d} \ | \ |p|<e^{a}(1-\delta)\}.
\end{equation}
Hence, the second term in the right side of (\ref{ghpln}) can be written
as
\begin{equation}
\label{ghpln2}
\frac{\widehat{G}(p)}{\hbox{ln}\Big(\frac{|p|}{e^{a}}\Big)}\chi_{A_{\delta}^{c +}}+
\frac{\widehat{G}(p)}{\hbox{ln}\Big(\frac{|p|}{e^{a}}\Big)}\chi_{A_{\delta}^{c -}}.
\end{equation}
Let us estimate using (\ref{inf1})
$$
\frac{|\widehat{G}(p)|}{|\hbox{ln}\Big(\frac{|p|}{e^{a}}\Big)|}\chi_{A_{\delta}^{c +}}
\leq \frac{\|G(x)\|_{L^{1}({\mathbb R}^{d})}}
{(2\pi)^{\frac{d}{2}}\hbox{ln}(1+\delta)}<\infty
$$
as we assume.

\noindent
Similarly,
$$
\frac{|\widehat{G}(p)|}{|\hbox{ln}\Big(\frac{|p|}{e^{a}}\Big)|}\chi_{A_{\delta}^{c -}}
\leq \frac{\|G(x)\|_{L^{1}({\mathbb R}^{d})}}
{-(2\pi)^{\frac{d}{2}}\hbox{ln}(1-\delta)}<\infty.
$$
We express
$$
\widehat{G}(p)=\widehat{G}(e^{a}, \sigma)+\int_{e^{a}}^{|p|}
\frac{\partial \widehat{G}(s, \sigma)}{\partial s}ds.
$$
Here and further down $\sigma$ will denote the angle variables on the sphere.
This allows us to write the first term in the right side of (\ref{ghpln}) as
\begin{equation}
\label{ghpln1}
\frac{\widehat{G}(e^{a}, \sigma)}{\hbox{ln}\Big(\frac{|p|}{e^{a}}\Big)}
\chi_{A_{\delta}}+\frac{\int_{e^{a}}^{|p|}
\frac{\partial \widehat{G}(s, \sigma)}{\partial s}ds}
{\hbox{ln}\Big(\frac{|p|}{e^{a}}\Big)}\chi_{A_{\delta}}.
\end{equation}
By virtue of the definition of the standard Fourier transform, we easily derive
\begin{equation}
\label{fdrub}
\Big|\frac{\partial \widehat{G}(p)}{\partial |p|}\Big|\leq
\frac{1}{(2\pi)^{\frac{d}{2}}}\|xG(x)\|_{L^{1}({\mathbb R}^{d})}, \quad p\in
{\mathbb R}^{d}, \quad d\geq 2.
\end{equation}
Let us estimate the second term in (\ref{ghpln1}) from above in the absolute
value by
$$
\frac{1}{(2\pi)^{\frac{d}{2}}}\|xG(x)\|_{L^{1}({\mathbb R}^{d})}
\Bigg|\frac{|p|-e^{a}}{\hbox{ln}\Big(\frac{|p|}{e^{a}}\Big)}\Bigg|\chi_{A_{\delta}}
\leq C\|xG(x)\|_{L^{1}({\mathbb R}^{d})}<\infty
$$
as assumed. Here and below $C$ will stand for a finite, positive constant.

\noindent
Thus, it remains to consider the term
\begin{equation}
\label{ghead}  
\frac{\widehat{G}(e^{a}, \sigma)}{\hbox{ln}\Big(\frac{|p|}{e^{a}}\Big)}
\chi_{A_{\delta}}.
\end{equation}
Clearly, (\ref{ghead}) is bounded if and only if $\widehat{G}(e^{a}, \sigma)$
vanishes, which is equivalent to orthogonality relations (\ref{or2}).

\noindent
We complete the proof of the lemma by considering the case a), such that the
dimension of the problem $d=1$. Let us introduce the intervals on the
real line
\begin{equation}
\label{idpm}
I_{\delta}^{+}:=[e^{a}(1-\delta),  e^{a}(1+\delta)], \quad
I_{\delta}^{-}:=[-e^{a}(1+\delta),  -e^{a}(1-\delta)],
\end{equation}
with $0<\delta<1$, such that
$$
I_{\delta}:=I_{\delta}^{+}\cup I_{\delta}^{-}.
$$
Obviously,
$$
I_{\delta}^{c}:=I_{\delta}^{c +}\cup I_{\delta}^{c -},
$$
where
\begin{equation}
\label{idcp}
I_{\delta}^{c +}:=(-\infty, -e^{a}(1+\delta))\cup (e^{a}(1+\delta), +\infty),  
\end{equation}
\begin{equation}
\label{idcm}
I_{\delta}^{c -}:=(-e^{a}(1-\delta), e^{a}(1-\delta)).
\end{equation}  
Evidently, we can write
\begin{equation}
\label{ghpln4}
\frac{\widehat{G}(p)}{\hbox{ln}\Big(\frac{|p|}{e^{a}}\Big)}=
\frac{\widehat{G}(p)}{\hbox{ln}\Big(\frac{|p|}{e^{a}}\Big)}\chi_{I_{\delta}^{+}}+
\frac{\widehat{G}(p)}{\hbox{ln}\Big(\frac{|p|}{e^{a}}\Big)}\chi_{I_{\delta}^{c +}}+
\frac{\widehat{G}(p)}{\hbox{ln}\Big(\frac{|p|}{e^{a}}\Big)}\chi_{I_{\delta}^{-}}+
\frac{\widehat{G}(p)}{\hbox{ln}\Big(\frac{|p|}{e^{a}}\Big)}\chi_{I_{\delta}^{c -}}.
\end{equation}
By means of (\ref{inf1}), we obtain the upper bound on the second term in the
right side of (\ref{ghpln4}) in the absolute value as
$$
\frac{|\widehat{G}(p)|}{\Big|\hbox{ln}\Big(\frac{|p|}{e^{a}}\Big)\Big|}
\chi_{I_{\delta}^{c +}}\leq \frac{\|G(x)\|_{L^{1}({\mathbb R})}}{\sqrt{2\pi}\hbox{ln}
(1+\delta)}<\infty
$$
due to our assumption.

\noindent
Analogously, for the fourth term in the right side of
(\ref{ghpln4}) we have the trivial inequality
$$
\frac{|\widehat{G}(p)|}{\Big|\hbox{ln}\Big(\frac{|p|}{e^{a}}\Big)\Big|}
\chi_{I_{\delta}^{c -}}\leq \frac{\|G(x)\|_{L^{1}({\mathbb R})}}{-\sqrt{2\pi}\hbox{ln}
(1-\delta)}<\infty.
$$
Note that
$$
\widehat{G}(p)=\widehat{G}(e^{a})+\int_{e^{a}}^{p}\frac{d\widehat{G}(s)}{ds}ds.
$$
Hence, the first term in the right side of (\ref{ghpln4}) equals to
\begin{equation}
\label{ghlni2}  
\frac{\widehat{G}(e^{a})}{\hbox{ln}\Big(\frac{|p|}{e^{a}}\Big)}
\chi_{I_{\delta}^{+}}+\frac
{\int_{e^{a}}^{p}\frac{d\widehat{G}(s)}{ds}ds}
{\hbox{ln}\Big(\frac{|p|}{e^{a}}\Big)}\chi_{I_{\delta}^{+}}.
\end{equation}
Using the definition of the standard Fourier transform (\ref{ft}), we easily
derive
\begin{equation}
\label{fdrub1}
\Big|\frac{d\widehat{G}(p)}{dp}\Big|\leq \frac{1}{\sqrt{2\pi}}\|xG(x)\|_
{L^{1}({\mathbb R})}.
\end{equation}
Clearly, formula (\ref{fdrub1}) is the analog of bound (\ref{fdrub}) in
one dimension.

\noindent
Let us obtain the estimate from above in the absolute value on the second term 
in (\ref{ghlni2}) as
$$
\frac{1}{\sqrt{2\pi}}\|xG(x)\|_{L^{1}({\mathbb R})}\Bigg|\frac{p-e^{a}}
{\hbox{ln}\Big(\frac{|p|}{e^{a}}\Big)}\Bigg|\chi_{I_{\delta}^{+}}\leq
C\|xG(x)\|_{L^{1}({\mathbb R})}<\infty     
$$
as we assume.

\noindent
Obviously, the first term in (\ref{ghlni2}) is bounded if and only if
$\widehat{G}(e^{a})=0$, which is equivalent to the orthogonality condition
$$
\Big(G(x), \frac{e^{ie^{a}x}}{\sqrt{2\pi}}\Big)_{L^{2}({\mathbb R})}=0.
$$
Evidently,
$$
\widehat{G}(p)=\widehat{G}(-e^{a})+\int_{-e^{a}}^{p}
\frac{d\widehat{G}(s)}{ds}ds.
$$
Thus, the third term in the right side of (\ref{ghpln4}) is equal to
\begin{equation}
\label{ghlni3}  
\frac{\widehat{G}(-e^{a})}{\hbox{ln}\Big(\frac{|p|}{e^{a}}\Big)}
\chi_{I_{\delta}^{-}}+\frac
{\int_{-e^{a}}^{p}\frac{d\widehat{G}(s)}{ds}ds}
{\hbox{ln}\Big(\frac{|p|}{e^{a}}\Big)}\chi_{I_{\delta}^{-}}.
\end{equation}
The second term in (\ref{ghlni3}) can be easily bounded from above in the
absolute value by
$$
\frac{1}{\sqrt{2\pi}}\|xG(x)\|_{L^{1}({\mathbb R})}\Bigg|\frac{p+e^{a}}
{\hbox{ln}\Big(\frac{-p}{e^{a}}\Big)}\Bigg|\chi_{I_{\delta}^{-}}\leq C
\|xG(x)\|_{L^{1}({\mathbb R})}<\infty     
$$
as assumed.

\noindent
The first term in (\ref{ghlni3}) is bounded if and only if 
$\widehat{G}(-e^{a})$ vanishes. This is equivalent to the orthogonality
relation
$$
\Big(G(x), \frac{e^{-ie^{a}x}}{\sqrt{2\pi}}\Big)_{L^{2}({\mathbb R})}=0,
$$
which completes the proof of the lemma. \hfill\lanbox

\bigskip

\noindent
We define the following technical expressions, which will help to
treat equations (\ref{id1m}).
\begin{equation}
\label{Nam}
N_{a,  \ m}:=
\Bigg\|\frac{\widehat{G_{m}}(p)}{\hbox{ln}\Big(\frac{|p|}{e^{a}}\Big)}\Bigg\|_
{L^{\infty}({\mathbb R}^{d})},
\end{equation}
where $a\in {\mathbb R}$ is a constant, $d\in {\mathbb N}$ and
$m\in {\mathbb N}$ as well.
Our final auxiliary proposition is as follows.

\bigskip

\noindent
{\bf Lemma A2.} {\it Let $m\in {\mathbb N}, \ G_{m}(x): {\mathbb R}^{d}\to
{\mathbb R}, \ G_{m}(x)\in L^{1}({\mathbb R}^{d}), \ xG_{m}(x)\in
L^{1}({\mathbb R}^{d})$ with $d\in {\mathbb N}$ and $a\in {\mathbb R}$, such
that 
$G_{m}(x)\to G(x)$ in $L^{1}({\mathbb R}^{d})$ and
$xG_{m}(x)\to xG(x)$ in $L^{1}({\mathbb R}^{d})$ as $m\to \infty$.

\medskip

\noindent
a) If $d=1$, let the orthogonality conditions
\begin{equation}
\label{or3}
\Big(G_{m}(x),\frac{e^{\pm ie^{a}x}}{\sqrt{2\pi}}\Big)_{L^{2}({\mathbb R})}=0, \quad m\in
{\mathbb N} 
\end{equation}
hold.

\medskip

\noindent
b) If $d\geq 2$, let the orthogonality relations
\begin{equation}
\label{or4}
\Big(G_{m}(x),\frac{e^{ipx}}{(2\pi)^{\frac{d}{2}}}\Big)_{L^{2}({\mathbb R}^{d})}=0,
\quad p\in S_{e^{a}}^{d}, \quad m\in
{\mathbb N} 
\end{equation}
be valid.

\medskip

\noindent
Let in addition
\begin{equation}
\label{2rnabml0}  
(2 \pi)^{\frac{d}{2}}N_{a, \ m}l\leq 1-\varepsilon
\end{equation}
for all $m\in {\mathbb N}$ with some fixed $0<\varepsilon<1$.

\medskip

\noindent
Then}
\begin{equation}
\label{Gmpaib}  
\frac{\widehat{G_{m}}(p)}{\hbox{ln}\Big(\frac{|p|}{e^{a}}\Big)} \to
\frac{\widehat{G}(p)}{\hbox{ln}\Big(\frac{|p|}{e^{a}}\Big)},
\quad m\to \infty ,
\end{equation}
{\it in $L^{\infty}({\mathbb R}^{d})$, so that}
\begin{equation}
\label{Gmpaibc}
\Bigg\|\frac{\widehat{G_{m}}(p)}{\hbox{ln}\Big(\frac{|p|}{e^{a}}\Big)}\Bigg\|_
{L^{\infty}({\mathbb R}^{d})}\to
\Bigg\|\frac{\widehat{G}(p)}{\hbox{ln}\Big(\frac{|p|}{e^{a}}\Big)}\Bigg\|_
{L^{\infty}({\mathbb R}^{d})}, \quad
m\to \infty.
\end{equation}
{\it Furthermore,
\begin{equation}
\label{Nabeps}  
(2 \pi)^{\frac{d}{2}}N_{a}l\leq 1-\varepsilon.
\end{equation}
}

\medskip

\noindent
{\it Proof.} By means of bound (\ref{inf1}), we have
\begin{equation}
\label{GmGp}
\|\widehat{G_{m}}(p)-\widehat{G}(p)\|_{L^{\infty}({\mathbb R}^{d})}\leq
\frac{1}{(2\pi)^{\frac{d}{2}}}\|G_{m}(x)-G(x)\|_{L^{1}({\mathbb R}^{d})}\to 0, \quad
m\to \infty
\end{equation}
as assumed.

\noindent
Evidently, (\ref{Gmpaibc})  follows easily
from (\ref{Gmpaib}) via the standard triangle inequality.

\noindent
Estimate (\ref{Nabeps}) is an easy consequence of (\ref{2rnabml0}) and  
(\ref{Gmpaibc}).

\noindent
Let us first consider the situation b) of our lemma when the dimension of
the problem $d\geq 2$.
In such case orthogonality conditions (\ref{or4}) hold as we assume. We easily
establish that the analogous statement will be valid in the limit. Obviously,
for $p\in S_{e^{a}}^{d}$
$$
\Big|\Big(G(x), \frac{e^{ipx}}{(2\pi)^{\frac{d}{2}}}\Big)_{L^{2}({\mathbb R}^{d})}\Big|=
\Big|\Big(G(x)-G_{m}(x), \frac{e^{ipx}}{(2\pi)^{\frac{d}{2}}}\Big)_{L^{2}({\mathbb R}^{d})}
\Big|\leq \frac{1}{(2\pi)^{\frac{d}{2}}}
\|G_{m}(x)-G(x)\|_{L^{1}({\mathbb R}^{d})}\to 0
$$
as $m\to \infty$ due to the one of our assumptions, so that
\begin{equation}
\label{G1l}
\Big(G(x), \frac{e^{ipx}}{(2\pi)^{\frac{d}{2}}}\Big)_{L^{2}({\mathbb R}^{d})}=0, \quad
p\in S_{e^{a}}^{d}
\end{equation}  
holds. We have
\begin{equation}
\label{gmhgpm}  
\frac{\widehat{G_{m}}(p)}{\hbox{ln}\Big(\frac{|p|}{e^{a}}\Big)}-
\frac{\widehat{G}(p)}{\hbox{ln}\Big(\frac{|p|}{e^{a}}\Big)}=
\frac{\widehat{G_{m}}(p)-\widehat{G}(p)}{\hbox{ln}\Big(\frac{|p|}{e^{a}}\Big)}
\chi_{A_{\delta}}+
\frac{\widehat{G_{m}}(p)-\widehat{G}(p)}{\hbox{ln}\Big(\frac{|p|}{e^{a}}\Big)}
\chi_{A_{\delta}^{c}}.
\end{equation}
Clearly, the second term in the right side of (\ref{gmhgpm}) equals to
\begin{equation}
\label{gmhgpmr2}
\frac{\widehat{G_{m}}(p)-\widehat{G}(p)}{\hbox{ln}\Big(\frac{|p|}{e^{a}}\Big)}
\chi_{A_{\delta}^{c +}}+
\frac{\widehat{G_{m}}(p)-\widehat{G}(p)}{\hbox{ln}\Big(\frac{|p|}{e^{a}}\Big)}
\chi_{A_{\delta}^{c -}}.
\end{equation}
Let us recall estimate (\ref{GmGp}). Hence
$$
\Bigg\|\frac{\widehat{G_{m}}(p)-\widehat{G}(p)}{\hbox{ln}\Big(\frac{|p|}{e^{a}}
\Big)}\chi_{A_{\delta}^{c +}}\Bigg\|_{L^{\infty}({\mathbb R}^{d})}\leq
\frac{\|G_{m}(x)-G(x)\|_{L^{1}({\mathbb R}^{d})}}
{(2\pi)^{\frac{d}{2}}\hbox{ln}(1+\delta)}\to 0, \quad  m\to\infty
$$
as assumed. Similarly,
$$
\Bigg\|\frac{\widehat{G_{m}}(p)-\widehat{G}(p)}{\hbox{ln}\Big(\frac{|p|}{e^{a}}
\Big)}\chi_{A_{\delta}^{c -}}\Bigg\|_{L^{\infty}({\mathbb R}^{d})}\leq
\frac{\|G_{m}(x)-G(x)\|_{L^{1}({\mathbb R}^{d})}}
{-(2\pi)^{\frac{d}{2}}\hbox{ln}(1-\delta)}\to 0, \quad  m\to\infty.
$$
By virtue of orthogonality relations (\ref{G1l}) and (\ref{or4}), we obtain
\begin{equation}
\label{gmgh0}  
\widehat{G}(e^{a}, \sigma)=0, \quad \widehat{G_{m}}(e^{a}, \sigma)=0, \quad
m\in {\mathbb N}.
\end{equation}
By means of (\ref{gmgh0}),
$$
\widehat{G}(p)=\int_{e^{a}}^{|p|}\frac{\partial \widehat{G}(s, \sigma)}
{\partial s}ds, \quad
\widehat{G_{m}}(p)=\int_{e^{a}}^{|p|}\frac{\partial \widehat{G_{m}}(s, \sigma)}
{\partial s}ds, \quad  m\in {\mathbb N}.     
$$
Then the first term in the right side of (\ref{gmhgpm}) is given by
\begin{equation}
\label{igmghln}
\frac{\int_{e^{a}}^{|p|}\Big[\frac{\partial \widehat{G_{m}}(s, \sigma)}{\partial s}-
\frac{\partial \widehat{G}(s, \sigma)}{\partial s}\Big]ds}
{\hbox{ln}\Big(\frac{|p|}{e^{a}}\Big)}\chi_{A_{\delta}}, \quad m\in {\mathbb N}.  
\end{equation}  
Let us use the definition of the standard Fourier transform (\ref{ft}) to
derive 
\begin{equation}
\label{dgmdgp}
\Big|\frac{\partial \widehat{G_{m}}(p)}{\partial |p|}-
\frac{\partial \widehat{G}(p)}{\partial |p|}\Big|\leq
\frac{1}{(2\pi)^{\frac{d}{2}}}\||x|G_{m}(x)-|x|G(x)\|_{L^{1}({\mathbb R}^{d})}, \quad
p\in {\mathbb R}^{d}, \quad d\geq 2, \quad m\in {\mathbb N},
\end{equation}
such that expression (\ref{igmghln}) can be easily bounded from above in the
absolute value by
$$
\frac{1}{(2\pi)^{\frac{d}{2}}}\||x|G_{m}(x)-|x|G(x)\|_{L^{1}({\mathbb R}^{d})}
\Bigg|\frac{|p|-e^{a}}{\hbox{ln}\Big(\frac{|p|}{e^{a}}\Big)}\Bigg|\chi_{A_{\delta}}
\leq C\||x|G_{m}(x)-|x|G(x)\|_{L^{1}({\mathbb R}^{d})}, \quad m\in {\mathbb N}.
$$
Therefore,
$$
\Bigg\|\frac{\int_{e^{a}}^{|p|}\Big[\frac{\partial \widehat{G_{m}}(s, \sigma)}
{\partial s}-
\frac{\partial \widehat{G}(s, \sigma)}{\partial s}\Big]ds}
{\hbox{ln}\Big(\frac{|p|}{e^{a}}\Big)}\chi_{A_{\delta}}\Bigg\|_{L^{\infty}({\mathbb R}^{d})}
\leq C\||x|G_{m}(x)-|x|G(x)\|_{L^{1}({\mathbb R}^{d})}\to 0, \quad m\to \infty         $$
as we assume. This completes the proof of the case b) of the lemma.

\noindent
Finally, we discuss the situation a) when the dimension of the problem $d=1$.
Hence, orthogonality relations (\ref{or3}) are valid as assumed. Let us
demonstrate that the analogous conditions will hold in the limit. Indeed,
$$
\Big|\Big(G(x), \frac{e^{\pm ie^{a}x}}{\sqrt{2\pi}}\Big)_{L^{2}({\mathbb R})}\Big|=
\Big|\Big(G(x)-G_{m}(x), \frac{e^{\pm ie^{a}x}}{\sqrt{2\pi}}\Big)_{L^{2}({\mathbb R})}
\Big|\leq \frac{1}{\sqrt{2\pi}}\|G_{m}(x)-G(x)\|_{L^{1}({\mathbb R})}\to 0
$$
as $m\to \infty$ due to the one of our conditions, such that
\begin{equation}
\label{or3l}  
\Big(G(x), \frac{e^{\pm ie^{a}x}}{\sqrt{2\pi}}\Big)_{L^{2}({\mathbb R})}=0.
\end{equation}
Obviously,
$$
\frac{\widehat{G_{m}}(p)}{\hbox{ln}\Big(\frac{|p|}{e^{a}}\Big)}-
\frac{\widehat{G}(p)}{\hbox{ln}\Big(\frac{|p|}{e^{a}}\Big)}=
\frac{\widehat{G_{m}}(p)-\widehat{G}(p)}{\hbox{ln}\Big(\frac{|p|}{e^{a}}\Big)}
\chi_{I_{\delta}^{+}}+
$$
\begin{equation}
\label{gmhpghpmln}
\frac{\widehat{G_{m}}(p)-\widehat{G}(p)}{\hbox{ln}\Big(\frac{|p|}{e^{a}}\Big)}
\chi_{I_{\delta}^{c +}}+
\frac{\widehat{G_{m}}(p)-\widehat{G}(p)}{\hbox{ln}\Big(\frac{|p|}{e^{a}}\Big)}
\chi_{I_{\delta}^{-}}+
\frac{\widehat{G_{m}}(p)-\widehat{G}(p)}{\hbox{ln}\Big(\frac{|p|}{e^{a}}\Big)}
\chi_{I_{\delta}^{c -}}, \quad m\in {\mathbb N}.
\end{equation}
By virtue of (\ref{GmGp}), the second term in the right side of
(\ref{gmhpghpmln}) can be estimated from above in the norm as
$$
\Bigg\|\frac{\widehat{G_{m}}(p)-\widehat{G}(p)}{\hbox{ln}\Big(\frac{|p|}
{e^{a}}\Big)}\chi_{I_{\delta}^{c +}}\Bigg\|_{L^{\infty}({\mathbb R})}\leq \frac
{\|G_{m}(x)-G(x)\|_{L^{1}({\mathbb R})}}{\sqrt{2\pi}\hbox{ln}(1+\delta)}\to 0, \quad
m\to \infty.            
$$
Similarly, the fourth term in the right side of
(\ref{gmhpghpmln}) can be bounded from above in the norm as
$$
\Bigg\|\frac{\widehat{G_{m}}(p)-\widehat{G}(p)}{\hbox{ln}\Big(\frac{|p|}
{e^{a}}\Big)}\chi_{I_{\delta}^{c -}}\Bigg\|_{L^{\infty}({\mathbb R})}\leq \frac
{\|G_{m}(x)-G(x)\|_{L^{1}({\mathbb R})}}{-\sqrt{2\pi}\hbox{ln}(1-\delta)}\to 0, \quad
m\to \infty.            
$$
Let us recall orthogonality relations (\ref{or3l}) and (\ref{or3}). Hence,
\begin{equation}
\label{gheam}  
\widehat{G}(e^{a})=0, \quad \widehat{G_{m}}(e^{a})=0, \quad m\in {\mathbb N}.
\end{equation}
By means of (\ref{gheam}), we have
$$
\widehat{G}(p)=\int_{e^{a}}^{p}\frac{d\widehat{G}(s)}{ds}ds, \quad
\widehat{G_{m}}(p)=\int_{e^{a}}^{p}\frac{d\widehat{G_{m}}(s)}{ds}ds, \quad
m\in {\mathbb N},
$$
so that the first term in the right side of (\ref{gmhpghpmln}) is equal to
\begin{equation}
\label{dgmhdhea}
\frac{\int_{e^{a}}^{p}\Big[\frac{d\widehat{G_{m}}(s)}{ds}-\frac{d\widehat{G}(s)}
{ds}\Big]ds}{\hbox{ln}\Big(\frac{|p|}{e^{a}}\Big)}\chi_{I_{\delta}^{+}}.  
\end{equation}  
Using the definition of the standard Fourier transform (\ref{ft}), we
easily arrive at
\begin{equation}
\label{dgmdgpabs}
\Bigg|\frac{d\widehat{G_{m}}(p)}{dp}-\frac{d\widehat{G}(p)}{dp}\Bigg|\leq
\frac{1}{\sqrt{2\pi}}\|xG_{m}(x)-xG(x)\|_{L^{1}({\mathbb R})}, \quad
p\in {\mathbb R}, \quad m\in {\mathbb N}.
\end{equation}  
Therefore, expression (\ref{dgmhdhea}) can be estimated from above in the
absolute value by 
$$
\frac{1}{\sqrt{2\pi}}\|xG_{m}(x)-xG(x)\|_{L^{1}({\mathbb R})}\Bigg|
\frac{p-e^{a}}{\hbox{ln}\Big(\frac{|p|}{e^{a}}\Big)}\Bigg|\chi_{I_{\delta}^{+}}\leq
C\|xG_{m}(x)-xG(x)\|_{L^{1}({\mathbb R})}, \quad m\in {\mathbb N},
$$
such that
$$
\Bigg\|\frac{\int_{e^{a}}^{p}\Big[\frac{d\widehat{G_{m}}(s)}{ds}-
\frac{d\widehat{G}(s)}{ds}\Big]ds}{\hbox{ln}\Big(\frac{|p|}{e^{a}}\Big)}
\chi_{I_{\delta}^{+}}\Bigg\|_{L^{\infty}({\mathbb R})}\leq C
\|xG_{m}(x)-xG(x)\|_{L^{1}({\mathbb R})}\to 0, \quad m\to \infty
$$
as assumed.

\noindent
By virtue of orthogonality conditions (\ref{or3l}) and (\ref{or3}), we have
\begin{equation}
\label{gheam-}  
\widehat{G}(-e^{a})=0, \quad \widehat{G_{m}}(-e^{a})=0, \quad m\in {\mathbb N}.
\end{equation}
Using (\ref{gheam-}), we easily obtain
$$
\widehat{G}(p)=\int_{-e^{a}}^{p}\frac{d\widehat{G}(s)}{ds}ds, \quad
\widehat{G_{m}}(p)=\int_{-e^{a}}^{p}\frac{d\widehat{G_{m}}(s)}{ds}ds, \quad
m\in {\mathbb N}.
$$
Then the third term in the right side of (\ref{gmhpghpmln}) is given by
\begin{equation}
\label{dgmhdhea-}
\frac{\int_{-e^{a}}^{p}\Big[\frac{d\widehat{G_{m}}(s)}{ds}-\frac{d\widehat{G}(s)}
{ds}\Big]ds}{\hbox{ln}\Big(\frac{|p|}{e^{a}}\Big)}\chi_{I_{\delta}^{-}}.  
\end{equation}
Let us use inequality (\ref{dgmdgpabs}) to derive the upper bound on
(\ref{dgmhdhea-}) in the absolute value by 
$$
\frac{1}{\sqrt{2\pi}}\|xG_{m}(x)-xG(x)\|_{L^{1}({\mathbb R})}\Bigg|
\frac{p+e^{a}}{\hbox{ln}\Big(\frac{|p|}{e^{a}}\Big)}\Bigg|\chi_{I_{\delta}^{-}}\leq
C\|xG_{m}(x)-xG(x)\|_{L^{1}({\mathbb R})}, \quad m\in {\mathbb N}.
$$
Then
$$
\Bigg\|\frac{\int_{-e^{a}}^{p}\Big[\frac{d\widehat{G_{m}}(s)}{ds}-
\frac{d\widehat{G}(s)}{ds}\Big]ds}{\hbox{ln}\Big(\frac{|p|}{e^{a}}\Big)}
\chi_{I_{\delta}^{-}}\Bigg\|_{L^{\infty}({\mathbb R})}\leq C
\|xG_{m}(x)-xG(x)\|_{L^{1}({\mathbb R})}\to 0, \quad m\to \infty
$$
due to the one of our assumptions, which completes the proof of the lemma in
the situation a).

\noindent
Note that under the given conditions
$$
N_{a}<\infty, \quad N_{a,  \ m}<\infty, \quad m\in {\mathbb N}, 
$$
where $a\in {\mathbb R}$ and $d\in {\mathbb N}$ according to
the result of Lemma A1 above. 
\hfill\lanbox 

\bigskip


\section*{Acknowledgements}

The first author is grateful to Israel Michael Sigal for the partial support
by the NSERC grant NA 7901.
The second author has been supported by the RUDN University Strategic Academic
Leadership Program.

\bigskip


\end{document}